\documentstyle[12pt]{article}
\begin{document}

\newcommand{\uA}{\mbox{$\underline{A}$}}
\newcommand{\uB}{\mbox{$\underline{B}$}}
\newcommand{\uS}{\mbox{${S}$}}
\newcommand{\de}{\mbox{$\delta $}}
\newcommand{\e}{\mbox{${\bf \epsilon \rm}$}}
\newcommand{\qed}{\mbox{$ \qquad \qquad \qquad \qquad \qquad  \qquad \qquad \qquad \qquad
\qquad \qquad \qquad \qquad \qquad \qquad \qquad \qquad
Q.E.D $ }}
\newtheorem{cor}{Corollary}[section]
\newtheorem{alg}{Algorithm}[section]
\newtheorem{lemma}{Lemma}[section]
\newtheorem{theo}{Theorem}[section]
\newtheorem{defi}{Definition}[section]
\newtheorem{exa}{Example} [section]
\newtheorem{pro}{Proposition}[section]
\newtheorem{rem}{Remark}[section]
\noindent{\bf OPEN-LOOP LINEARIZATION OF NON-LINEAR DISCRETE INPUT - OUTPUT SYSTEMS THROUGH SIMPLIFICATION  ALGORITHMS. \rm}
\vskip 10 pt
\noindent
{\it Stelios Kotsios, \par\noindent  University of Athens, Department of Economics, 
\par\noindent
Division of Mathematics and Computer Science
\par\noindent
Pesmazoglou 8, 10559 Athens, Greece}
\vskip 10 pt
\par\noindent
\small
{\bf Abstract:} The problem of linear equivalence for a general class of nonlinear systems, is examined throughout this paper. A relevant algorithm is developed, based on a factorization procedure. This factorization is based on the star-product, an operation corresponding to the cascade connection of systems.
\vskip 5 pt
\noindent
{\bf Keywords:} Nonlinear, Discrete, Factorization, Algebraic, Computational, Linearization, Simplification.
\normalsize
\section{Introduction}
\noindent
As it is known, nonlinear systems are used in a variety of applications and have been the focus of research for a number of years. Many of mathematical tools have been used for the analysis and design of such systems. Among these tools, algebraic and computational methods have had a major effect. These tools have led to the development of new algorithms and techniques, which in turn have allowed design methods to be accomplished with greater speed and efficiency, \cite{kn:carcanias}.
\par
\noindent
Such computational methods have been applied in the study of modeling,
the problem of feedback linearization, and in the global optimization problem.
(\cite{kn:glad3},\cite{kn:fliessnew},\cite{kn:carcanias}, to mention but a few).
 All of these problems were investigated mostly in relation to continuous systems, and, in some cases, in relation to non-linear discrete
 time systems, \cite{kn:factn}.
By definition, a discrete time system evaluates input and output signals over a
countable number of time instants. In some cases, a discrete time system is obtained
from continuous time systems through sampling at certain time instants
\cite{kn:gladbook}. A further  case concerns  systems that are naturally and directly
described in discrete form, typically in financial or economic systems, \cite{kn:puu}.
There is a rich literature devoted to the study of discrete systems. 
Certain works
approach the issue through analytical tools, \cite{kn:cyrot1},\cite{kn:monaco} 
and others by using algebraic methodologies, like
differential algebra or rings theory, 
\cite{kn:glad3},\cite{kn:fliessnew},
\cite{kn:glad1},\cite{kn:perdon},\cite{kn:wu}.
\par
\noindent
Factorization is a popular algebraic method dealing with linear as well as nonlinear, continuous or discrete, systems. In linear systems the central idea is that of factorizing the transfer matrix as the ratio of two stable rational matrices, \cite{kn:gladbook}. In nonlinear systems construction of coprime factorizations have been studied by several authors. This setting was mainly applied to systems that can be described by input-output representations,\cite{kn:fliessnew},\cite{kn:sontag}.
\par
\noindent
The above factorizations can be used effectively in many design problems, \cite{kn:sontag},\cite{kn:hammer}. The simplification-linearization method is among them. This method can permit us to determine systems with a " less " complex structure, which are equivalent to the original nonlinear system, i.e. they will give the same output under the same input and identical initial conditions. Among these " simpler " systems the linear systems are the most desirable ones, because their behavior is well known. Therefore, we may ask when a nonlinear system is equivalent to a linear one. This is a version of the linearization problem.
\par
\noindent
\par\noindent
The present paper focuses on the simplification-linearization problem of nonlinear discrete input-output systems of the form:
\[  y(t)+ \sum a_i y(t-i)+\sum \sum a_{i_1i_2} y(t-i_1) y(t-i_2) + \cdots +\]
\[+\sum
\sum \cdots \sum a_{i_1i_2 \ldots i_n} y(t-i_1) \cdots y(t-i_n)=\]
\[=\sum b_j u(t-j) +\cdots +
\sum \sum \cdots \sum b_{j_1j_2 \ldots j_m} u(t-j_1) \cdots u(t-j_m)+ \]
\[+\sum\sum c_{kl}y(t-k)u(t-l)+\cdots +\sum \sum\sum \cdots \sum\sum
\cdots \sum \]
\begin{equation}\label{main}
                 c_{k_1k_2 \ldots k_{n'}l_1 l_2 \ldots
l_v} y(t-{k_1}) \cdots y(t-{k_{n'}}) u(t-l_1) \cdots u(t-l_v)
\end{equation}
Equation (\ref{main}) transforms causal input signals (i.e. $u(t)=0$ for $t<0$) to causal
output signals. A set of initial conditions
$y(0)=y_0$, $y(1)=y_1$, $\ldots$, $y(k)=y_k$ is always assigned to (\ref{main}) and the lowest
delay output term (that is $y(t)$) appears in the linear part of the system.
The systems of the form (\ref{main}), which contain products among inputs and output signals that are
sometimes called 'cross-products', encompass a broad variety of nonlinear discrete
systems. We obtain these either through transformations of nonlinear discrete
state-space representations into input - output forms \cite{kn:rugh}, or when we use
Taylor's expansion method to approximate other more general nonlinear discrete systems
\cite{kn:kalouptsidis},\cite{kn:rugh}. These are employed in signal processing theory,
whenever it is necessary to construct nonlinear representations of discrete signals,
(they are an extension of the infinite impulse response filters to a nonlinear set up, 
\cite{kn:kalouptsidis}),
in nonlinear time-series analysis and in adaptive control, in the context of  designing
nonlinear adaptive controllers, \cite{kn:kokotovic}.
\par
\noindent
In order to describe the above systems algebraically we create an appropriate framework by using the so-called $\de\e$-operators to deal with cross-products.
This
operator, introduced in \cite{kn:bibo}, is an extension of the simple $\de$ and
$\e$-operators that are used in the algebraic description of nonlinear discrete input
- output systems without cross-products \cite{kn:phd1}, and which are also an
extension of the simple shift operator $q$ in linear discrete systems,
\cite{kn:amstrong}. Via these operators we can define the so-called $\de$,$\e$ and $\de \e$ - polynomials. By means of these polynomials we can formally rewrite (\ref{main})
as follows: $ Ay(t)=Bu(t)+{C}[y(t),u(t)]$, where $A$ a $\de$-polynomial, $B$
an $\e$-polynomial and ${C}$ a $\de \e$-polynomial. Among these polynomials, we can define two
product operations: the dot-product (denoted by '$\cdot$'), which corresponds to the usual product among
polynomials, and the star-product (denoted by '$\ast$'), which corresponds to the substitution of one
polynomial by another, or, in more system-oriented terminology, to the cascade
connection of systems. Though there is a similarity between this algebraic background
and others, there is also the interesting peculiarity that the set of $\de
\e$-polynomials with respect to the star-product is not a ring \cite{kn:glad3}.
\par
\noindent
Two systems $A_1y(t)=B_1u(t)+{C}_1[y(t),u(t)]$ and
$A_2y^*(t)=B_2v(t)+{C}_2[y^*(t),v(t)]$ are equivalent, if $y(t)=y^*(t)$,
whenever $u(t)=v(t)$ and initial conditions are identical. A " simplification "
 method for a given non-linear system consists in discovering systems that are equivalent to the original one but have a " less " complex
 structure. This procedure allows us to replace, if necessary, the original system with the simpler one. Clearly, of all such 'simpler' systems, linear systems are the most desirable ones. Accordingly, in this paper we develop a method 
 for discovering linear systems which are equivalent to a given non-linear system.
In other words, we are looking for a linear system $A_ly^*(t)=B_lv(t)$, such that
$y^*(t)=y(t)$, whenever $u(t)=v(t)$. In other words we have a kind of the linearization procedure devoted to open-loop systems, thus, in our approach we do not use any feedback law design.
\par
\noindent
To face the simplification-linearization problem we separate the systems of the form (\ref{main}) into two categories, the first consists from nonlinear systems without cross-products, the second from nonlinear systems with cross-products. 
To work with systems of the first class we
follow the next steps: First, be means of a symbolic algorithm, we factorize the nonlinear system as $\hat{A} \ast L y(t)=\hat{B} \ast Mu(t)$, where
$L,M$ are linear $\de$ and $\e$-polynomials correspondingly and $\hat{A},\hat{B}$ nonlinear. The novelty is that the coefficients of all the above polynomials are not constant number but functions of certain undetermined parameters $w_{ij},s_{ij}$. By giving to these parameters certain values, we endow the polynomials with specific properties. Therefore, using certain methods of computational algebra we can find values of $w_{ij},s_{ij}$ so that $\hat{A}=\hat{B}$ and hence the nonlinear polynomial can be simplified to the linear one $Ly(t)=Mu(t)$. We have to mention here that this 
simplification procedure is something which is well known in the
linear case. Indeed, let $q$ be the simple shift operator, i.e.
$q^{-1}y(t)=y(t-1)$ and $Ay(t)=Bu(t)$ a linear system, $A,B$ polynomials
of the variable $q^{-1}$. If we re-write
it in the form
$C \cdot A'y(t)=C \cdot B'u(t)$ and $C$ is a stable polynomial, then, by performing a division,
we take the linear system
$ A'y(t)=B'u(t)$, which gives the same output with the original system, under the same input and
initial conditions \cite{kn:amstrong}. 
To deal with polynomials of the second class, that is those with cross-products, we work in an analogous way: First we factorize the system as $\underline{\hat{A}}\ast [L,M]=0$, where $L,M$ linear $\de$ and $\e$-polynomials and $\underline{\hat{A}}$ a $\de\e$-polynomial with parametrical coefficients. Proper choice of values for the parameters may give to $\underline{\hat{A}}$ a certain structure, the so-called homogeneous structure. If this happens, then we can write the nonlinear system as $\tilde{A}\ast (Ly(t)+Mu(t))$, $\tilde{A}$ a nonlinear polynomial, and then by simplification of $\tilde{A}$ we take the linear equivalent system $Ly(t)=-Mu(t)$.
\par
\noindent
The particularities of our methodology are: 1) Its symbolic algorithmic orientation that allows a direct computer implementation. Actually, we create a proper algebraic framework so that specific algorithms to be applied. These algorithms permit us to transform questions from system theory, such as linearization, to computational algebra questions. Indeed, in order to take solutions to our simplification procedure we have to solve systems of polynomial equations, this can be done by means of certain algorithms of computational algebra. We note here, that our approach is symbolical and not numerical.
2) It gives not a single solution but a class of them. Indeed, if the above mentioned systems of polynomial equations do not accept a single solution but many of them, then to its such solution corresponds a linear system, equivalent to the nonlinear one.
\par
\noindent
Finally, we have to make clear, that in the present paper we are not involved with internal stability questions. We are mainly interest to describe the simplification procedures formally. 
The text is divided into two parts. In Part I we present the algebraic
framework and the algorithms, in Part II we are dealing separately, with the simplification-linearization 
problem of systems without cross-products and of systems which include cross-products. Throughout
the paper, ${\bf N,Z,R}$ will denote the set of natural, rational and real numbers,
respectively.

\section{Part I - The Algebraic Preliminaries}
In this section we present the algebraic preliminaries, needed for the development of the simplification algorithms.

\subsection{The Multiindices}

Let $k$ be a positive integer. A subset of the set $\cup^k_{n=1}Z^n$ is called a set
of {\bf multiindices} and it is denoted by ${\bf I}$. A set of multiindices may be finite or
not. We denote the elements of ${\bf I}$ by ${\bf i}=(i_1,i_2,\ldots,i_n)$. Usually, we put the 
elements of a multiindex ${\bf i}$ in an ascending way, that is
$i_1 \le i_2 \le \cdots \le i_n $. The quantity
$i_1+i_2+ \cdots +i_n$ , called the degree of ${\bf i}$, is denoted by $deg({\bf i})$.
We can define the following operations among multiindices. These will be used in the
results in systems simplification below.
\begin{itemize}
\item
Let ${\bf i}=(i_1,i_2,\ldots,i_n)$ be an multiindex. An multiindex of the form ${\bf
j}=(i_{\sigma_1},i_{\sigma_2}, \ldots, i_{\sigma_k})$, with $i_{\sigma_\theta} \in \{
i_1,i_2,\ldots,i_n\}$ and $\sigma_1 \ne \sigma_2 \ne \cdots \ne \sigma_k $, is called
a $k$-subindex of ${\bf i}$.

\item
Let ${\bf i}=(i_1,i_2,\ldots,i_n)$ be an multindex and ${\bf
j}=(i_{\varphi_1},i_{\varphi_2}, \ldots,i_{\varphi_k})$ a subindex of ${\bf i}$. Their
difference is defined as follows: ${\bf i \ominus j}=(i_{\lambda_1},i_{\lambda_2}, \ldots,
i_{\lambda_m})$, with $i_{\lambda_\mu} $ $\in \{ i_1,i_2, \ldots, i_n\}$ $-\{
i_{\varphi_1},i_{\varphi_2}, \ldots, i_{\varphi_k} \}$ .

\item
Given two multiindices $ {\bf i \rm}=(i_1,i_2,\ldots,i_k) $ and $ {\bf j
\rm}=(j_1,j_2,\ldots,j_{\lambda}) $, the new multiindex $ {\bf i \oplus j \rm} $ is defined
just juxtaposing $ {\bf j \rm  }$ after $ {\bf i \rm} $. Explicitly, $ {\bf i \oplus j
\rm}=(i_1,j_1,i_2,i_3,j_2,
 \ldots,i_k,j_{\lambda}) $, where $  i_1 \le j_1 \le i_2 \le
i_3 \le j_2 \le \ldots \le i_k \le j_{\lambda}$.

\item
Let $ {\bf j \rm}=(j_1,j_2,\ldots,j_m) $ be a multiindex. We define the pointwise sum $
{\bf j \rm} \dot + \varphi $, $\varphi$ a real number, as follows: $ {\bf j \rm}  \dot
+ \varphi=(j_1+\varphi,j_2+\varphi,\ldots,j_m+\varphi) $.
\end{itemize}
Finally, a set ${\bf I \rm}$ of multiindices can be ordered in a lexicographical
way as follows: Let $ {\bf i \rm}= (i_1,i_2, \ldots, i_n ) $ 
and $ \bf j \rm =(j_1,j_2,\ldots,j_m)$ be two multiindices, with elements ordered 
in an increasing way. We can say that the multiindex $ \bf i \rm  $  is "less" than the
multiindex $ \bf j $, and we write $ {\bf i < j \rm } $, if either
 $ n<m  $ or
$ n=m $ and the right-most nonzero entry of the vector ${\bf j} - {\bf i}$ is 
positive.
Since for any two multiindices we can
immediately deduce which multiindex is "larger" and which "smaller", the above
lexicographical order is well defined.

\subsection{The $\de\e$-operators}
We introduce the $\de \e$-operator to deal with cross-products among input and
output signals, which acts upon a pair of sequences. This initially appeared in
\cite{kn:bibo} and has been examined in \cite{kn:lappas1},\cite{kn:lappas2}. Let
$y(t),u(t)$ be real sequences, defined over the set of integers ${\bf Z}$  and let $F$ be
the set of causal sequences. Thus if $y(t),u(t) \in F $ then $y(t)=u(t)=0$
 for $ t<0$.
Let ${\bf i}=(i_1, i_2, \ldots, i_m), {\bf j }=(j_1, j_2, \ldots, j_n) $ be multiindices.
The operator $ \de_{\bf i} \times \de_{\bf j} : F \times F \mapsto F $ is defined as:
$ \de_{\bf i} \times \de_{\bf j}{[y(t), u(t)]}=$ $ y(t-i_1)y(t-i_2) \cdots
y(t-i_m)u(t-j_1)u(t-j_2) \cdots u(t-j_n) $. This means that the operator $\de_{\bf i}$
acts exclusively on " outputs " and $\de_{\bf j}$ exclusively on " inputs ".
Sometimes, the following  notation may be more convenient: $ \de_{\bf i} \times
\de_{\bf j}=\de_{\bf i} \epsilon_{\bf j} $. We call these operators $ \de \e
$-operators and we denote their set by $\Delta \times \Delta$, where by $\Delta$ we
denote the set of simple $\de_{\bf i}$ operators. By convention, we define  $ \delta_e
\{ y(t) \}=\{1\} $ for each $ t \in Z $. This means that $ \de_{\bf i} \times \de_e
[y(t), u(t)]=\de_{\bf i}y(t) $, $ \de_e \times \de_{\bf j} [y(t),u(t)]=\e_{\bf j}u(t) $.
Practically speaking, the symbol $e$ corresponds to the empty multiindex $\{ \}$.
Obviously, the operators $\de_i \e_e=\de_i$ or $\de_e \e_i=\e_i$, $i \in Z$, are
 none other than the well-known
simple shift operators, i.e., 
$\de_i \e_e y(t)=y(t-i)$ or $\de_e \e_i u(t)=u(t-i) $.
The operators $\de_0,\e_0$, give
zero delays, i.e. $\de_0\e_ey(t)=y(t)$, $\de_e\e_0u(t)=u(t)$. The $\de\e$-operators of the
form $\de_0\de_{i_2}\cdots \de_{i_m}\e_0\e_{j_2}\cdots \e_{j_n}$ are many times called zero terms.
\begin{exa}
Let ${\bf i}=(0,1,1,2)$, ${\bf j}=(1,1,3,3,3)$, ${\bf r}=(0,0,1,2)$ and ${\bf
s}=(0,1,1,2)$, then $ \de_{\bf i} \e_{\bf j}[y(t),u(t)]=$ $\de_{(0,1,1,2)}$
$\e_{(1,1,3,3,3)}[y(t),u(t)]=$ $\de_0\de_1\de_1\de_2$ $\e_1\e_1\e_3\e_3\e_3
[y(t),u(t)]$ $=\de_0\de_1^2\de_2$ $\e_1^2\e_3^3[y(t),u(t)]=$
$y(t)y^2(t-1)y(t-2)u^2(t-1)u^3(t-3)$, furthermore $\de_{\bf
r}\e_{e}[y(t),u(t)]=\de_{\bf r}y(t)=$ $\de_0^2\de_1\de_2y(t)=y^2(t)y(t-1)y(t-2)$ and
$\de_e \e_{\bf s}[y(t),u(t)]=\e_{\bf s}u(t)=$ $\e_0\e_1^2\e_2 u(t)=u(t)u^2(t-1)u(t-2)$
\end{exa}
The order among multiindices implies an order among $\de \e$-operators in a natural way.
Indeed, we say that $\de_{\bf i}\e_{\bf j} \preceq \de_{\bf i'}\e_{\bf j'}$ if either
${\bf i} < {\bf i'}$ or ${\bf i} = {\bf i'}$ and ${\bf j} < {\bf j'}$. We equip the
set of the $\de \e$-operators with two internal operations: the dot-product and the
star-product. The dot-product corresponds to the usual product among sequences, while
the star-product corresponds to the substitution of one sequence by another.
Specifically, if $z(t)=\de_{\bf i}\e_{\bf j} [y(t),u(t)]$, $s(t)=\de_{\bf i'}\e_{\bf
j'} [y(t),u(t)]$, ${\bf i}=(i_1,i_2, \ldots, i_n)$, ${\bf j}=(j_1,j_2, \ldots, j_m)$,
${\bf i'}=(i'_1,i'_2, \ldots, i'_{n'})$, ${\bf j'}=(j'_1,j'_2, \ldots, j'_{m'})$, then
$ \de_{\bf i}\e_{\bf j}$ $ \cdot \de_{\bf i'}\e_{\bf j'}[y(t),u(t)]=$ $y(t-i_1) \cdot
y(t-i_2)$ $ \cdots y(t-i_n)\cdot y(t-i'_1) \cdots y(t-i'_{n'})$ $ u(t-j_1) \cdot
u(t-j_2) $ $\cdot u(t-j_m)\cdot u(t-j'_1) \cdots u(t-j'_{m'})$. Let us now suppose that
$z(t)=\de_{\bf i}\e_{\bf j} [y(t),u(t)]$, and $y(t)=\de_{\bf i'}w(t)$, $u(t)=\e_{\bf
j'}v(t)$,
 then
$ \de_{\bf i}\e_{\bf j} \ast [\de_{\bf i'}, \e_{\bf j'}][w(t),v(t)]=$ $\de_{\bf
i}\e_{\bf j}[\de_{\bf i'}w(t), \e_{\bf j'}v(t)]$. The meaning of the last
relation is the following:
Let us define the maps
 $Z: F \times F \to F,$$Z[y(t),u(t)]=z(t)=\de_{\bf i}\e_{\bf j}[y(t),u(t)]$,
 $Y: F \to F, Y[w(t)]=y(t)=\de_{\bf i'}w(t) $,
 $U:F \to F, U[v(t)]=u(t)=\e_{\bf j}v(t)$. 
 The star-product is nothing else than the composition
$Z \circ [Y,U]=Z[Y,U]$ or by using the other notation
 $\de_{\bf i}\e_{\bf j} \ast [\de_{\bf i'},\e_{\bf j'}]=$
$\de_{\bf i}\e_{\bf j} \circ [\de_{\bf i'},\e_{\bf j'}]$$=\de_{\bf i}\e_{\bf j}[\de_{\bf i'},
\e_{\bf j'}]$.

\par
\noindent The following properties are valid for all the above operations. Their
proofs are obvious from the definitions and are therefore omitted.
\begin{pro}\label{proti}
(a) $ \de_{\bf i}\e_{\bf j} \cdot \de_{\bf i'}\e_{\bf j'}=$ $\de_{\bf i
\oplus i'} \e_{\bf j \oplus j'}$ (b) $ \de_{\bf i}\e_{\bf j} \ast [\de_{\bf i'},
\e_{\bf j'}]=$ $\de_{{\bf i'} \dot{+}i_1} \de_{{\bf i'} \dot{+}i_2} \cdots \de_{{\bf
i'} \dot{+}i_n} $ $ \e_{{\bf j'} \dot{+}j_1} \e_{{\bf j'} \dot{+}j_2} \cdots  \e_{{\bf
j'} \dot{+}j_m}$ (c) $\de_{\bf i}\e_{\bf j} \ast ([y_1(t),u_1(t)]+[y_2(t),u_2(t)])$ $
\neq \de_{\bf i}\e_{\bf j} \ast [y_1(t),u_1(t)] + \de_{\bf i}\e_{\bf j} \ast
[y_2(t),u_2(t)]$
\end{pro}
Part (c) of the above proposition, indicates that the $\de\e$-operators do not behave
like members of an algebraic ring.
\begin{exa} Let ${\bf i}=(0,1)$, ${\bf j}=(1,1,2)$, ${\bf i'}=(0,0,1)$ and ${\bf j'}=(0,0)$.
Then, $\de_{\bf i} \e_{\bf j}$ $ \cdot \de_{\bf i'} \e_{\bf j'}[y(t),u(t)]=$
$\de_0\de_1 \e_1^2\e_2$ $ \de_0^2\de_1\e_0^2 [y(t),u(t)]=$
$\de_0^3\de_1^2\e_0^2\e_1^2\e_2 [y(t),u(t)]$ $=y^3(t)y(t-1)u^2(t)u^3(t-1)u(t-2)$.
Since ${\bf i} \oplus {\bf i'}=(0,0,0,1,1)$ and $ {\bf j} \oplus {\bf j'}=(0,0,1,1,2)$
we can easily verify formula (a) of proposition (\ref{proti}). Moreover, $\de_{\bf
i}\e_{\bf j} \ast [\de_{\bf i'}, \e_{\bf j'}][y(t),u(t)]=$ $\de_0\de_1\e_1^2\e_2 \ast
[y^2(t)y(t-1),u^2(t)]$ $=y^2(t)y(t-1)y^2(t-1)y(t-2)u^4(t-1)u^2(t-2)$. Part (b), of
proposition (\ref{proti}), can be verified too, through straight substitution.
\end{exa}

\subsection{The $\de \e$-polynomials}
The $\de \e$-polynomials are straightforward extensions of the $\de \e$-operators that
are necessary for the description of discrete non-linear systems with products among
input and output sequences. These have been introduced and studied in
\cite{kn:bibo},\cite{kn:lappas1}. We will go over the main points, adding some new
results that will clarify the ideas in the present paper. Expressions of the form:
$A=\sum_{({\bf i,j}) \in {\bf I \times J}} c_{\bf ij}\de_{\bf i} \e_{\bf j}$, are
called $\de \e$-polynomials, where by ${\bf I, J}$ we denote sets of multiindices and
$c_{\bf ij}$ real numbers. By putting ${\bf I}=\{\de_e\}$, or ${\bf J}=\{ \e_e \}$, we can
transform the above polynomials to pure $\de$ or $\e$-polynomials. That is
$A=\sum_{\bf i \in I} c_{{\bf i}e}\de_{\bf i}$ or $A=\sum_{\bf j \in J}
c_{e{\bf j}}\e_{\bf j}$. These are polynomials that deal only with input or output
signals and produce non-linear polynomial expressions of these, without cross-products.
By $\de_{\bf i} A$ if $A$ is a $\de$-polynomial, or $\e_{\bf j} A$ if $A$ is
an $\e$-polynomial, we denote the following operations:
$ \de_{\bf i} A=A_{i_1}\cdot A_{i_2}\cdots A_{i_m}$ 
with $ A_{i_k}=\sum_{\bf i \in I} c_{{\bf i}e}\de_{{\bf i}\dot{+}i_k}$
and
$ \e_{\bf j} A=A_{j_1}\cdot A_{j_2}\cdots A_{j_n}$
with $ A_{j_k}=\sum_{\bf j \in J} c_{e{\bf j}}\e_{{\bf j}\dot{+}j_k}$. In other words we shift all the delays of $A$ by $i_1$, then by $i_2$ and so on and finally we multiply the obtained polynomials.
We will call expressions of the form $\sum_{i=1}^k c_{ ie}\de_{ i}$ or $\sum_{
j=1}^{\lambda} c_{e j}\e_{ j}$, linear $\de$ or $\e$-polynomials, respectively. The
term, which accordingly to the order defined previously is ordered highly among the
terms of $A$, is called the {\it maximum} term of $A$. By $d(A)$ we denote the minimum
delay of $A$, in other words
 $d(A)=\min\{\min({\bf \rho }), {\bf \rho} \in {\bf I} \cup {\bf J}\} $.
The largest of the numbers
 $deg({\bf i}+{\bf j})$, is called the {\it degree} of $A$, denoted also by $deg(A)$. 
 Two $\de \e$-polynomials $A$ and $A'$ are equal if
$A[y(t),u(t)]=A'[y(t),u(t)]$ for each $y(t),u(t) \in F$ and $t\ge 0$. If $A$ is a $\de
\e$-polynomial and $B$, $C$ are a $\de$ and an $\e$-polynomial, then their
star-product, $A \ast [B,C]$, is defined as the composition $A \circ [B,C]$. This
means that $B$ is replaced into the $\de$-part of $A$ and $C$ into the $\e$-part. Formulae
for the calculation of the star-product and additional properties, can be found in
\cite{kn:bibo}, \cite{kn:phd1}, \cite{kn:lappas1}. These papers prove that the set of
the pure $\de$ or $\de \e$-polynomials does not form an algebraic ring. This property
makes them different to all other similar approaches \cite{kn:glad3}, \cite{kn:glad1}.
The following formulas will be used in the proofs of the main theorems of this paper,
below. Their proofs are straightforward and are therefore omitted.
\begin{pro}\label{star}
(a) Let $A$ be a $\de$-polynomial and $B$ an $\e$-one, then 
$ \de_{\bf i} \e_{\bf j} \ast [A,B]=(\de_{\bf i} A) 
\cdot (\e_{\bf j} B)$
\par
\noindent (b) Let $\de_{\bf i}\e_{\bf j}$ be a $\de$-operator, with ${\bf
i}=(i_1,i_2,\ldots,i_{\varphi})$, ${\bf j}=(j_1,j_2,\ldots,j_{\lambda})$ and
$L=\sum_{i \in I_1}r_i \de_i $, $M=\sum_{j \in J_1} m_j \e_j$ linear $\de$ and
$\e$-polynomials, $I_1,J_1 \subset Z$ finite sets of multiindices, then
\[ \de_{\bf i}\e_{\bf j} \ast [L,M]=\sum_{{\bf s} \in I_1^{\varphi}}
\sum_{{\bf h} \in J_1^{\lambda}}r_{s_1}r_{s_2}\cdots r_{s_\varphi} m_{h_1} m_{h_2}
\cdots
 m_{h_\lambda}\de_{{\bf s}{+}{\bf i}}\e_{{\bf h}{+}{\bf j}} \]
 where ${\bf s}=(s_1,s_2,\ldots, s_{\varphi})$, ${\bf h}=(h_1,h_2, \ldots, h_{\lambda})$ and
 $I_1^{\varphi}, J_1^{\lambda}$ the cartesian products
 \newline
 $I_1^{\varphi}=\underbrace{I_1\times I_1 \times \cdots \times I_1}_{\varphi - times}$,
 $J_1^{\lambda}=\underbrace{J_1\times J_1 \times \cdots \times J_1}_{\lambda - times}$.
 \end{pro}
 Property (a) is a remarkable one. It creates, by suitable choosing of 
 $n$ and $m$, any polynomial term including products of certain delays
 among input and output signals. Therefore, it indicates that the star-product
 acts likewise the derivative in the case of continuous systems, (it creates terms
 with derivatives of any order), \cite{kn:kotritt}.
  \begin{exa}
 (1) $ \de_1\de_2 \e_2^2 \ast [2\de_0+\de_1,2\e_0-\e_1]=$
 $[\de_1\de_2 \ast (2\de_0+\de_1) ] \cdot [\e_2^2\ast(2\e_0-\e_1)]=$
 $(2\de_1+\de_2) \cdot (2\de_2+\de_3) $
 $\cdot (2\e_2-\e_3) \cdot (2\e_2-\e_3)=$
 $16\e_2^2\de_1\de_2 - 16\e_2\e_3\de_1\de_2 +$
$ 4\e_3^2\de_1\de_2 + 8\e_2^2\de_2^2 - 8\e_2\e_3\de_2^2 +$ $2\e_3^2\de_2^2 +
8\e_2^2\de_1\de_3 - $ $ 8\e_2\e_3\de_1\de_3 + 2\e_3^2\de_1\de_3 +$ $ 4\e_2^2\de_2\de_3
- 4\e_2\e_3\de_2\de_3 + \e_3^2\de_2\de_3$
\par
(2) $(2 \de_1\e_2+\de_0^2\e_1\e_2)$ $ \ast [\de_1-2\de_2, \e_0^2\e_1]$ $=2 \de_1\e_2
\ast [\de_1-2\de_2, \e_0^2\e_1]$ $+ \de_0^2\e_1\e_2 \ast [\de_1-2\de_2, \e_0^2\e_1]$
$=[2 \de_1 \ast (\de_1-2\de_2)] \cdot [\e_2 \ast \e_0^2\e_1]$ $+[\de_0^2 \ast
(\de_1-2\de_2)] \cdot [\e_1\e_2 \ast  \e_0^2\e_1]$ $=2(\de_2-2\de_3) \cdot
(\e_2^2\e_3)+$ $(\de_1-2\de_2)^2 \cdot \e_1^2\e_2 \cdot \e_2^2\e_3=$
$\e_1^2\e_2^3\e_3\de_1^2 + $ $2\e_2^2\e_3\de_2 +$ $ 4\e_1^2\e_2^3\e_3\de_2^2 -$
$4\e_2^2\e_3\de_3$$-4\de_1\de_2\e_1^2\e_2^3\e_3$

\end{exa}
The so-called "homogeneous " polynomials are special cases of $\de \e$-polynomials with specific
applications.
\begin{defi}
Let ${\bf \theta}=(i_1,i_2, \ldots, i_n)$ be a given multiindex. We call a $\de
\e$-polynomial $T$, "homogeneous" with respect to the multiindex ${\bf \theta}$ and the constants
$\lambda, \mu$, if it has the form: $T=\mu^n \de_e\e_{i_1}\e_{i_2} \cdots \e_{i_n} + $
$ \lambda \mu^{n-1} \de_{i_1}$ $ \e_{i_2} \e_{i_3} \cdots \e_{i_n}$ $+\lambda
\mu^{n-1} \de_{i_2}$ $ \e_{i_1} \e_{i_3} \cdots \e_{i_n}+ \cdots $ $+\lambda \mu^{n-1}
\de_{i_n}$ $ \e_{i_1} \e_{i_2} \cdots \e_{i_{n-1}}+$ $\lambda^2 \mu^{n-2}
\de_{i_1}\de_{i_2}$ $ \e_{i_3}\e_{i_4} \cdots \e_{i_n}$ $+\lambda^2 \mu^{n-2}
\de_{i_1}\de_{i_3}$ $ \e_{i_2}\e_{i_4} \cdots \e_{i_n}$ $+\cdots +\lambda^2 \mu^{n-2}$
$ \de_{i_{n-1}}\de_{i_n} \e_{i_1}\e_{i_2}$ $ \cdots \e_{i_{n-2}}+$ $\cdots
+\lambda^{n-1} \mu \de_{i_1}\de_{i_2}$ $ \cdots \de_{i_{n-1}} \e_{i_n} $ $+ \cdots
+\lambda^{n-1} \mu \de_{i_2}\de_{i_3}$ $ \cdots \de_{i_{n}} \e_{i_1}+$ $\lambda^n
\de_{i_1} \de_{i_2} \cdots \de_{i_n}\e_e =$
\begin{equation}\label{round}
=\sum_{\varphi =0}^n\sum_{\bf (i,j), i \in I_{\varphi}, j=\theta \ominus i}
\lambda^{\varphi}\mu^{n-\varphi} \de_{\bf i} \e_{\bf j}
\end{equation}
where $\lambda, \mu \in R$ and ${\bf I}_{\varphi}$ is the set of the
$\varphi$-subindices of the multiindex ${\bf i}$, ${\bf I}_0=\{ e \}$.
\end{defi}
The following property is of special interest.
\begin{pro}\label{roundpro}
Let $A$ be a $\de \e$-polynomial. If $A$ can be written in the form: $A=\tilde{A} \ast
[L,M]$, where $L,M$ are linear $\de$ and $\e$-polynomials respectively and $\tilde{A}$
is a homogeneous $\de \e$-polynomial, with respect to an multiindex ${\bf \theta }$ and the
constants $\lambda,\mu$, then $ A=\de_{\bf \theta} \ast (\lambda L +\mu M) $.
\end{pro}
{\bf Proof:} Using the $\de$-symbol, for the delays of $L$ and the $\e$-symbol, for
the delays of $M$, we get the following succession: $ \de_{\bf \theta} \ast (\lambda L
+\mu M)=$ $ \de_{i_1} \de_{i_2} \cdots \de_{i_n} \ast (\lambda L +\mu M)=$ $ \de_{i_1}
\ast (\lambda L +\mu M) \cdot $ $ \de_{i_2} \ast (\lambda L +\mu M)$ 
$\cdots$ $ \de_{i_n} \ast (\lambda L +\mu M)=$ $\lambda^n (\de_{i_1}
\ast L) \cdot (\de_{i_2} \ast L)$ $\cdots (\de_{i_n} \ast L)+\cdots +$ $\lambda
\mu^{n-1} (\de_{i_1} \ast L) \cdot (\e_{i_2} \ast M) \cdots $ $(\e_{i_n} \ast M ) +
\cdots $ $+ \mu^n (\e_{i_1} \ast M ) \cdot (\e_{i_2} \ast M ) \cdots $ $(\e_{i_n} \ast
M ) $ $=\cdots=$
$\lambda^n \de_{i_1} \de_{i_2} \cdots \de_{i_n}\e_e \ast [L,M]+ \cdots+\Large$ $\lambda
\mu^{n-1} \de_{i_1}\e_{i_2} \cdots \e_{i_n} \ast [L,M] +\cdots +$ $\mu^n \de_e \e_{i_1}
\e_{i_2} \cdots \e_{i_n} \ast [L,M]=$ $[\lambda^n \de_{i_1} \de_{i_2} \cdots
\de_{i_n}\e_e + \cdots +$ $\lambda \mu^{n-1} \de_{i_1}\e_{i_2} \cdots \e_{i_n}  +\cdots +$
$\mu^n \de_e \e_{i_1} \e_{i_2} \cdots \e_{i_n}] \ast [L,M]=$ $\tilde{A} \ast [L,M]$,
where $\tilde{A}$ is a homogeneous polynomial with respect to the multiindex ${\bf \theta}$ and
the constants $\lambda, \mu$. The proposition has been proved. $\Box$

\subsection{The Formal $\de \e L$-Factorization.}
In this section we shall present the main tool that is used throughout this paper. Let
$\{w_{ij}\}$ and $\{ s_{ij} \}$, $i=1,2, \ldots $, $j=1,2, \ldots $ be two sets of
parameters, i.e.,  sets of variables that can take arbitrary real values. Let $A$ be a $\de
\epsilon$-polynomial. We call 'Formal-$\de \epsilon$-Linear-Like-Factorization' of $A$
($F\de \e L$-Factorization for short),
 any expression of the form:
\[ A=\de_0 \e_e \ast [A_{\de ,
l},\e_e]+\de_e \e_0 \ast [\de_e, A_{\e , l}]+\sum_{k=1}^g c_k(w_{ij},s_{ij})\de_{{\bf i}_k}
\e_{{\bf j}_k} \ast [L_k,M_k] +R\] 
where $c_k(w_{ij},s_{ij})$ are functions of the
parameters $w_{ij},s_{ij}$, 
$A_{\de , l},
A_{\e , l}$ the $\de$ and $\e$-linear terms of $A$,
$\de_{{\bf i}_k}\e_{{\bf j}_k} \in \Delta \times \Delta$,
with $\min({\bf i}_k)=0$, $\min({\bf j}_k)=0$ and 
\[ L_k=w_{k0}\de_0+w_{k1}\de_1+\cdots+\de_k \quad,\quad M_k=s_{k0}\e_0+s_{k1}\de_1+\cdots+\e_k \]
That is, they are linear $\de$ and
$\e$-polynomials, where their coefficients are not constant numbers but 
the parameters $w_{ij},s_{ij}$ correspondingly (by the exception of the highest delay term). The $\de \epsilon$-polynomial $R$, with
parametrical coefficients,
 which is called
the {\it remainder}, must contain only {\it zero} terms, i.e., terms of the form $\de_0
\de_{i_1} \ldots \de_{i_{\phi}}$ $\epsilon_0 \epsilon_{j_1} \ldots
\epsilon_{j_{\sigma}}$. The $F \de \e L$-Factorization of a given polynomial $A$ is
denoted by $Formal[A]$. For pure $\de$-polynomials or
$\epsilon$-polynomials, we can restrict the above terminology by either
$\epsilon_{j_k}=\epsilon_e$, $M_k=\epsilon_e$ or $\de_{j_k}=\de_e$, $L_k=\de_e$. For
instance, we say that a $\de$-polynomial $B$, has a Formal - $\de$-Linear - Like -
Factorization
 or, simply, an $F\de L$-Factorization,
 if $B=\sum_{k=1}^g c_k(w_{ij})\de_{{\bf i}_k}\e_e \ast [L_k,\e_e] +R $ or in short
$ B=\sum_{k=1}^g c_k\de_{{\bf i}_k} \ast L_k +R $. (See \cite{kn:difeq} for further
details.)
\begin{exa}\label{formal}
A formal $\de \e L $-Factorization of the $\de \e$-polynomial $B=\de_1^2+\de_1\e_1$ is
$ Formal[B]=\de_0^2 \ast (w_{10}\de_0$ $+\de_1) +\de_0 \e_0 \ast[w_{20}\de_0$ $+\de_1,
s_{20}\e_0+\e_1]-$ $s_{20} \de_0\e_0 \ast [w_{30}\de_0+\de_1,\e_0]$ $-w_{20}\de_0\e_0
\ast[\de_0,s_{40} \e_0+\e_1]+$ $[-2w_{10}\de_0\de_1-w_{10}^2 \de_0^2 +$
$(w_{20}s_{20}+w_{30}s_{20}$ $+w_{20}s_{40})\de_0\e_0]$, where
$w_{10},w_{20},s_{20},w_{30},s_{40}$ are the parameters.
\end{exa}

\begin{theo}
For a given $\de \e$-polynomial $A$, the $F \de \e L$-Factorization is unique, assuming
that the parameters $\{ w_{ij},s_{ij} \}$ are considered as constants.
\end{theo}
{\bf Proof:} We explain first that by the expression " are considered as 
constants " we mean that we treat the parameters like to be specific number.
Now, if $A$ has only zero terms then the proof of the theorem is trivial and
$c_{k}=0$ and $R=A$. Let us suppose that $A$ contains non-zero terms and
$\omega_a\de_{{\bf i}_a} \e_{{\bf j}_a}$ is its maximum term. Using proposition
(\ref{star}), (b), we can see that this term also appears in the product
$c_a(w_{ij},s_{ij})\de_{{\bf i}_a} \e_{{\bf j}_a} \ast [L_a,M_a]$
$=c_a(w_{ij},s_{ij})\de_{{\bf i}_a} \e_{{\bf j}_a} \ast$ $[w_{0a}\de_0+w_{1a}\de_1+
\cdots +\de_{m_a}, s_{0a}\e_0+s_{1a}\e_1+\cdots+\e_{n_a}]$. By equating their
coefficients we calculate the quantity $c_a(w_{ij},s_{ij})$ uniquely, in effect
$c_a(w_{ij},s_{ij})=\omega_a$. If we repeat this procedure for the next higher order
term, we find an expression for the "next" coefficient $c_{a-1}(w_{ij},s_{ij})$. Since
this expression is a function of $c_a(w_{ij},s_{ij})$ and of some parameters
$w_{ij},s_{ij}$ considered as constants, we conclude that $c_{a-1}(w_{ij},s_{ij})$ is
also defined uniquely. By induction, we can then see that all the coefficients are
uniquely determined. The polynomial $R$ consists only of zero terms. These terms arise
either from the polynomial $A$ or from the products $c_k(w_{ij},s_{ij})\de_{{\bf i}_k}
\e_{{\bf j}_k} \ast [L_k,M_k]$.
 The unique determination of the coefficients
 ${c}_{k}$ entails the uniqueness of $R$ and, thus,  the theorem has been proved.
$\Box$

\subsection{ The Algorithms.}
In what follows we provide the basic algorithmic tool that we will then use to obtain
a Formal $\de \e L$-Factorization. This is along the lines of Ritt's remainder
algorithm \cite{kn:ritt}. It is based on a kind of division with respect to the
star-product \cite{kn:kotritt}. \small {\sf
\par
\vskip 5 pt \underline{\bf The $F \de\e L  $-Subroutine.}
\par
\vskip 5 pt {\bf Input:} A non-linear $\de \e$-polynomial $A$, without linear or
constant terms.
\par
{\bf Initial Condition:} The index $\lambda=0$.
\par
\noindent \vskip 5 pt {\bf DO}
\par
\begin{quote}
{\bf STEP 1:} We set $\lambda=\lambda+1$.
\par
{\bf STEP 2:} Let $S=c_{\bf \lambda}(w_{ij},s_{ij})\de_{\bf i}\e_{\bf j}=c_{\bf
\lambda}(w_{ij},s_{ij}) \de_{i_1}\de_{i_2} \cdots \de_{i_n} \e_{j_1}\e_{j_2} \cdots
\e_{j_m} $ be the {\it maximum} non-zero term of $A$, $0 \le i_1 \le i_2 \le  \cdots
\le i_n$, $0 \le j_1 \le j_2 \le  \cdots \le j_m$ and $c_{\bf \lambda}(w_{ij},s_{ij})$
its coefficient. (Actually, in the first iteration of the algorithm, $c_{\bf
\lambda}(w_{ij},s_{ij})$ is always a real number. It then becomes a function of the
unknown parameters $w_{ij},s_{ij}$.)
\par
{\bf STEP 3:} We form the linear $\de$-polynomial
\[ L_\lambda=w_{\lambda 0}\de_0+w_{\lambda 1}\de_1+
w_{\lambda 2}\de_2 +\cdots +\de_{i_1} \]
where $w_{\lambda 0},w_{\lambda 1},\ldots$ 
are unknown parameters, taking values in
${\bf R}$.
\par
{\bf STEP 4:} We form the linear $\e$-polynomial
\[ M_\lambda=s_{\lambda 0}\e_0+s_{\lambda 1}\e_1+s_{\lambda 2}\e_2 
+\cdots +\e_{j_1} \]
where $s_{\lambda 0},s_{\lambda 1},\ldots$ are unknown parameters, taking values in
${\bf R}$.
\par
{\bf STEP 5:} We calculate the quantity:
\[ R_\lambda=A-c_{\bf \lambda}(w_{ij},s_{ij})\de_{0}\de_{i_2-i_1} \cdots \de_{i_n-i_1}
\e_{0}\e_{j_2-j_1} \cdots \e_{j_m-j_1}
 \ast [L_\lambda, M_\lambda]=\]
 \[=
A-c_{\bf \lambda}(w_{ij},s_{ij}) \de_{\bf i_{\lambda}} \e_{\bf j_{\lambda}}\ast
[L_\lambda, M_\lambda] \] {\bf STEP 6:} We replace the polynomial $A$ with the
polynomial $R_\lambda$.
\end{quote}
{\bf UNTIL} All the terms of $A$ become zero terms.
\par
\noindent {\bf STEP 7:} We rename the last value of $A$ as $R$.
\par
\noindent {\bf Output:} The quantities $c_{k}(w_{ij},s_{ij}),\de_{\bf i_{\rm
k}}\e_{\bf j_{\rm k}}, L_k, M_k, R$, $k=1,\ldots,\varphi$. \vskip 10 pt } {\sf
\par
\vskip 5 pt \noindent \underline{\bf The $F \de\e L  $-Algorithm.}
\par
\vskip 5 pt \noindent {\bf Input:} A $\de \e$-polynomial $A$.
\par
\vskip 5 pt \noindent
\begin{quote}
{\bf STEP 1:} We decompose $A$ as follows: $ A=A_{\de , l}+A_{\e , l}+A_{nl}$, where $
A_{\de , l} , A_{\e , l}, A_{nl}$ are the $\de$-linear, $\e$-linear and the non-linear
parts of $A$, respectively.
\par
\noindent {\bf STEP 2:} {\bf IF} $A_{nl}=0$ {\bf THEN} we set $\varphi = 0$, $R=0$
{\bf ELSE} using the $F \de \e L$-Subroutine we take the quantities
$c_{k}(w_{ij},s_{ij}),\de_{\bf i_{\rm k}}\e_{\bf j_{\rm k}}, L_k, M_k, R$,
$k=1,\ldots,\varphi$.
\par
\noindent {\bf STEP 3:} We set $c_{\varphi+1}(w_{ij},s_{ij})=1$, $\de_{{\bf
i}_{\varphi +1}}\e_{{\bf j}_{\varphi +1}}=\de_0\e_e$, $L_{\varphi+1}=A_{\de,l}$,
$M_{\varphi +1}=\e_e$, $c_{\varphi+2}(w_{ij},s_{ij})=1$, $\de_{{\bf i}_{\varphi
+2}}\e_{{\bf j}_{\varphi +2}}=\de_e\e_0$, $L_{\varphi+2}=\de_e$, $M_{\varphi
+1}=A_{\e,l}$.
\par
\end{quote}
\vskip 5 pt \noindent {\bf Output: }The quantities $c_{k}(w_{ij},s_{ij}),\de_{\bf
i_{\rm k}}\e_{\bf j_{\rm k}}, L_k, M_k, R$, $k=1,\ldots,\varphi+2$. \vskip 5 pt }

\normalsize

\begin{theo}
The $F\de \e L$-Algorithm terminates after a finite number of steps. If
$c_{k},\de_{\bf i_{\rm k}}\e_{\bf j_{\rm k}}, L_k, M_k, R$, $k=1,\ldots,\varphi+2$,
are its outputs, then the quantity $Formal[A]=\sum_{k=1}^{\varphi+2} c_k\de_{{\bf
i}_k}\e_{{\bf j}_k} \ast [L_k,M_k] +R$ is the Formal - $\de \e L$ - Factorization of
the given $\de \e$-polynomial $A$.
\end{theo}
{\bf Proof:} If we follow the above algorithm step by step, we get $
R_1=A-c_{1}(w_{ij},s_{ij})\de_{{\bf i}_1} \e_{{\bf j}_1}\ast [L_1,M_1] $. By means of
proposition (\ref{star}), (b), we conclude that this operation eliminates
 at least one non-zero term of $A$. Next, we get
 $ R_2=R_1-c_{2}(w_{ij},s_{ij})\de_{{\bf i}_2}\de_{{\bf j}_2} \ast [L_2,M_2] $,
 which eliminates another non-zero term of $A$ and so on. These repeated operations,
ensure that all the non-zero terms of $A$ will be finally eliminated and therefore the
procedure will terminate. This also means that the remainder $R$ will contain only
zero terms. Furthermore, for the operators $\de_{{\bf i}_\lambda}\e_{{\bf j}_\lambda}$
in Step 5 of the $F\de\e L$-Subroutine, we have $\min({{\bf i}_\lambda})=\min({{\bf
j}_\lambda})=0$. Additionally, by reverse substitution we get $
A=c_{1}(w_{ij},s_{ij})\de_{{\bf i}_1} \e_{{\bf j}_1}\ast [L_1,M_1]
 + c_{2}(w_{ij},s_{ij})\de_{{\bf i}_2}\de_{{\bf j}_2} \ast [L_2,M_2]
  +\cdots + R $
$=\sum_{k=1}^{\varphi +2} c_k(w_{ij},s_{ij})\de_{{\bf i}_k}\e_{{\bf j}_k} \ast
[L_k,M_k] +R$. Finally, all the above ensure that the output polynomials $\de_{\bf
i_{\rm k}}\e_{\bf j_{\rm k}}, L_k, M_k, R$, $k=1,\ldots,\varphi +2$ and the
coefficients $c_{k}(w_{ij},s_{ij})$, constitute a formal $\de \e L$-Factorization for
the given $\de \e $-polynomial $A$. $\Box$
\par\noindent

\subsection{Evaluations}
Let $A$ be a $\de \e$-polynomial and $Formal[A]$, its Formal $\de \e L$-Factorization.
If the parameters $w_{ij},s_{ij}$ take values according to a set of substitution rules
$U$, we say that the $Formal[A]$ is evaluated over the set $U$, thus writing: $\left.
\begin{array} {c}
Formal[A]\\
\end{array} \right|_{U} $.
More rigorously, let ${\bf W}=( w_{ij})$ be the set of the $w$-parameters, written as
a vector and ${\bf S}=( s_{kl})$ the vector of the $s$-parameters. ${\bf r}=( a_{ij}
)$, ${\bf q}=( b_{kl} )$ are vectors of real numbers, which are in one-to-one
correspondence with the vectors ${\bf W,S}$. We say that these parameters follow the
rule $({\bf r,q })$, thus writing ${\bf (W , S) \to (r,q)}$, if the following
substitutions are valid $ (w_{ij},s_{kl})=(a_{ij}, b_{kl})$. Let $N,Q$ be two sets of
rules, $N=\{ {\bf r}_1,{\bf r}_2, \ldots, {\bf r}_{\lambda}\}$, $Q=\{ {\bf q}_1,{\bf
q}_2, \ldots, {\bf q}_{\mu}\}$ and $U \subseteq N \times Q$, then
\[ \left. \begin{array}
{c}
Formal[A]\\
\end{array} \right|_{U} =\bigcup_{\sigma=1}^{\xi} \left\{ \sum_{k=1}^g c_k(w_{ij},s_{hl})
\de_{{\bf i}_k}\e_{{\bf j}_k} \ast [L_k,M_k] +R,\right.\]
\[  \quad with, \quad \left. {\bf (W,S) \to (r, q)_{\sigma} } \in U \right\} \]
The set of substitution rules, $U$, may be finite or infinite. Often we denote this by
$U=\{(w_{ij},s_{kl})=(a_{ij},b_{kl})_\sigma, \sigma = 1,2, \ldots \}$, or $U=\{
(w_{ij}=a_{ij},s_{kl}=b_{kl}),$ $(w_{ij}=a'_{ij},s_{kl}=b'_{kl}),\ldots,\}$. 
\begin{exa}
Let $B$ and $Formal[B]$ be as in example (\ref{formal}). The vectors of the parameters
are ${\bf W}=(w_{10},w_{20},w_{30})$, ${\bf S}=(s_{20},s_{40})$. For instance, let us
take the following sets of rules $N=\{ (1,-1,2)\}$, $Q=\{(-1, \frac{1}{2}) \}$, then
$U=\{ (w_{10},w_{20},w_{30}, s_{20},s_{40})$ $=(1,-1,2,-1, \frac{1}{2}) \}$ and $
\left. \begin{array} {c}
Formal[B]\\
\end{array} \right|_{U} =\de_0^2\ast (\de_0+\de_1)+\de_0\e_0 \ast [-\de_0+\de_1, -\e_0+\e_1]-$
$(-1)\de_0\e_0 \ast [2\de_0+\de_1,\e_0]$ $-(-1)\de_0\e_0 \ast [\de_0, \frac{1}{2} \e_0
+\e_1]$ $+( -2\de_0\de_1-\de_0^2-\frac{1}{2}\de_0\e_0)$. The rules $N=\{
(0,\varphi,-\varphi (\theta +1))\}$, $Q=\{ (1,\theta)\}$, will give $U=\{
(w_{10},w_{20},w_{30}, s_{20},s_{40})$ $=(0,\varphi,-\varphi (\theta +1), 1,\theta) $,
$ \varphi , \theta \in {\bf R}$ and $ \left. \begin{array} {c}
Formal[B]\\
\end{array} \right|_{U} =\de_0^2\ast \de_1$
$+\de_0\e_0 \ast [ \varphi \de_0+\de_1, \e_0+\e_1]-$ $ \de_0\e_0 \ast [ -\varphi
(\theta +1) \de_0+\de_1,\e_0]$ $- \varphi \de_0\e_0 \ast [\de_0,
 \theta \e_0 +\e_1]$.
 
 \end{exa}

\section{Part II-The Simplification-Linearization Algorithms}
In this section we present the simplification procedures of our methodology.
Before we examine them, 
we would like to begin by a short algebraic
description of the non-linear discrete input-output systems, 
via the notion of $\de
\e$-polynomials.

\subsection{The $\de \e$-polynomials and Non-linear Discrete Systems}

Suppose that we have a system of the form (\ref{main}). By using
$\de\e$-operators we can rewrite this as follows: 
\[\sum a_i \de_i y(t) + \sum \sum
a_{ij} \de_i \de_j y(t) + \cdots + \sum \cdots \sum a_{i_1 i_2 \ldots i_n}\de_{i_1}\de_{i_2} \cdots \de_{i_n}y(t)=\]
\begin{equation}\label{demain}
=\sum b_i \e_i u(t) +\sum \sum
b_{ij} \e_i \e_j u(t)+ \cdots+\sum \cdots \sum \de_{i_1} \cdots \de_{i_m}
\e_{j_1} \cdots \e_{j_k} [y(t),u(t)]
\end{equation} 
We can describe the above relation shortly, by writing
$Ay(t)=Bu(t)+{C}[y(t),u(t)]$, 
where $A$ is a $\de$-polynomial, $B$ an
$\e$-polynomial and ${C}$ a pure $\de \e$ - one. The causality and
solvability of the system is guaranteed by the inequality: $d(A_l)=d(A)< \min \{ d(B),
d({C}) \}$, where $A_l$ is the linear part of $A$.
This inequality means that the lower delayed term of the output (i.e. $y(t)$) appears
in the linear part of the system and thus we can solve (\ref{main}) with respect
to $y(t)$ in a direct way.

To each nonlinear system of the form (\ref{main}) we assign a vector
of real numbers ${\bf y}_0=(y_0,y_1, \ldots, y_{k-1})$ which
gives the so called initial conditions:
$y(0)=y_0$, $y(1)=y_1$, \ldots, $y(k-1)=y_{k-1}$, where $k$ is the maximum
delay appeared in the output signal. Since the signals, involved in
(\ref{main}), are causal, i.e. $y(t)=0$, $u(t)=0$, for $t<0$, for each
given vector of initial conditions any input signal $u(t)$ determines
 a unique output signal that satisfies (\ref{main}), for $n \ge k$.
Let us have two causal nonlinear systems: $Ay(t)=Bu(t)+{C}[y(t),u(t)]$
with ${\bf y}_0=(y_0,y_1, \ldots, y_{k-1})$ and
$\hat{A}\psi(t)=\hat{B}v(t)+{\hat{C}}[\psi(t),v(t)]$ with 
${\bf \psi}_0=(\psi_0,\psi_1,\ldots,\psi_{\lambda-1})$ and $k>\lambda$.
We say that the two systems " operate "
under the same initial conditions or that their initial conditions are identical
if $\psi_0=y_0$,
$\psi_1=y_1$, $\ldots$, $\psi_{\lambda-1}=y_{\lambda-1}$ 
and $y_{\lambda}=\psi(\lambda)$, 
$y_{\lambda+1}=\psi(\lambda+1)$, $\ldots$, 
$y_{k-1}=\psi(k-1)$. In other words we must give
 as initial conditions to the system which starts to product outputs
 later, the corresponding outputs of the other system.
\par
\noindent
Two systems $A_1y(t)=B_1u(t)+{C}_1[y(t),u(t)]$ and
$A_2y^*(t)=B_2v(t)+{C}_2[y^*(t),v(t)]$ are equivalent, if $y(t)=y^*(t)$,
whenever $u(t)=v(t)$ and initial conditions are identical. A " simplification "
 method for a given non-linear system consists in discovering systems that are equivalent to the original one but have a " less " complex
 structure. This procedure allows us to replace, if necessary, the original system with the simpler one.
\par
\noindent
A $\de\e$-polynomial is called 'proper' if the following two facts obtain: i) The
minimum delay of $A$, $\alpha=d(A)$, appears only in the $\de$-part of certain terms
of $A$. ii) The power of $\de_\alpha$, in the above terms, is equal to one. Obviously,
all polynomials with  linear parts that contain the lowest 
delay terms are proper. 
The
following theorem is an extension of a useful result that appears in \cite{kn:factn}.
The proof is along similar lines.
\begin{theo}\label{proper}
The equality $A[y(t),u(t)]=A[y^*(t),u(t)]$, $\forall t \in Z$, implies that
$y(t)=y^*(t)$, provided that $A$ is proper and $y(t),y^*(t)$ have the same initial
conditions.
\end{theo}
A direct consequence of the above theorem is the following result:
\begin{cor}\label{porisma}
Let $A$ be a proper $\de$-polynomial, then the equality $Ay(t)=0$, $\forall t \in Z$ implies that $y(t)=0$ under zero initial conditions.
\end{cor}

All the above terminology can be applied to systems without cross-products, in a direct way. 
These are systems of the form:
\[\sum a_i y(t-i) + \sum\sum
a_{ij} y(t-i)y(t-j) + \cdots + \sum \cdots \sum a_{i_1 i_2 \ldots i_n}y(t-{i_1})y(t-{i_2})\cdots y(t-{i_n})=\]
\[=\sum b_i u(t-i) +\sum \sum
b_{ij} u(t-i)u(t-j)+ \cdots+\sum \cdots \sum u(t-{j_1}) \cdots u(t-{j_k})\]
Using $\de$ and $\e$-operators we rewrite the above system as
\[\left(\sum a_i \de_i  + \sum\sum
a_{ij} \de_i \de_j  + \cdots + \sum \cdots \sum a_{i_1 i_2 \ldots i_n}\de_{i_1}\de_{i_2} \cdots \de_{i_n}\right)y(t)=\]
\[=\left(\sum b_i \e_i  +\sum \sum
b_{ij} \e_i \e_j + \cdots+\sum \cdots \sum \e_{j_1} \cdots \e_{j_k}\right) u(t)\]
or shortly $Ay(t)=Bu(t)$, where $A$ a $\de$-polynomial and $B$ an $\e$-polynomial. The causality is guaranteed
by the inequality: $d(A_l)=d(A)<d(B)$, where $A_l$ is the linear part of $A$.
Two systems $A_1y(t)=B_1u(t)$ and
$A_2y^*(t)=B_2v(t)$ are equivalent, if $y(t)=y^*(t)$,
whenever $u(t)=v(t)$ and initial conditions are identical. A " simplification "
procedure for these systems  is defined by the same way, as before.

\subsection{The Case without Cross-Products}
In this section we are dealing with systems which do not contain products among input and output signals. 
Actually, we want to simplify them via the use of a symbolic algorithm.
Clearly, of all such 'simpler' systems, linear systems are the most desirable ones. Accordingly, in this section we ask under what circumstances  a non-linear system of the form $Ay(t)=Bu(t)$ is equivalent to a linear one.
In other words, we are looking for a linear system $A_ly^*(t)=B_lv(t)$, such that
$y^*(t)=y(t)$, whenever $u(t)=v(t)$. 
We are not involved with stability questions.
The
algorithm upon discussion is: 
 
\small
\vskip 15 pt
{\sf
\underline{\bf The Linear-Equivalence Algorithm I}
\vskip 15 pt
{\bf Input:} The nonlinear $\de$, $\e$-polynomials $A$ and $B$.
\par
\begin{quote}
\par
{\bf STEP 1:} By means of the $F\de \e L $-Algorithm, we find the quantities:
\[ Formal[A]=\sum_{k=1}^\omega c_{k}(w_{ij}) \de_{{\bf i}_k} \ast L_k +R_{\de} \]
\[ Formal[B]=\sum_{\rho=1}^\varphi c_{\rho}'(s_{ij}) \e_{{\bf j}_\rho} \ast M_\rho +R_{\e} \]

\par
{\bf STEP 2:} We form the set of rules:
 ${U}=\{ (w_{ij},s_{ij})=(a_{ij}, b_{ij}), a_{ij},b_{ij} \in {\bf R} \}$, such that the following are valid simultaneously
 \begin{itemize}
  \item $\left. \begin{array}
{c}
R_{\de}\\
\end{array} \right|_{U}
= \left. \begin{array}
{c}
R_{\e}\\
\end{array} \right|_{U}=0$

\item The sets ${\bf L}=\{L: L=\gcd(\left. \begin{array}
{c}
L_k\\
\end{array} \right|_{U} ) \}$ 
and 
${\bf M}=\{M: M=\gcd(\left. \begin{array}
{c}
M_\rho\\
\end{array} \right|_{U} ) \}$ 
are not void, i.e.  ${\bf L} \ne \emptyset$ and ${\bf M} \ne \emptyset$.
\item The polynomials $\hat{A}, \hat{B}$, where 
$A=\hat{A} \ast L$, $L \in {\bf L}$, $B=\hat{B} \ast M, M \in {\bf M}$, are proper and
$\hat{A}=\hat{B}$.

  \end{itemize}
\par
{\bf STEP 3:} {\bf IF} ${U} \ne \emptyset$ {\bf THEN} 
goto the output {\bf ELSE} the method fails.
\end{quote}
\vskip 15 pt
{\bf Output:} The quantities ${\bf L,M}$

\normalsize

\begin{theo}\label{lin1}
Let us suppose that we have the causal nonlinear discrete system 
$Ay(t)=Bu(t)$, $A$, $B$ $\de$ and $\e$ - polynomials 
with $d(A)<d(B)$.
If ${\bf L}$ and ${\bf M}$, are the outputs of the previous algorithm, then each linear system of the form
$Ly(t)=Mu(t)$, $L \in {\bf L}$, $M \in {\bf M}$ is equivalent to
the original nonlinear system.
\end{theo}
{\bf Proof:} Let, ${\bf W}$ and ${\bf S}$, be the sets of parameters, 
appearing in the formal factorization of $A$, 
and $B$, during the run of the algorithm. 
The 
existence of the solution
implies that ${U} \ne \emptyset$
and therefore we can find at least, one pair of vectors, ${\bf r}=( a_{ij})$,
${\bf q}=( b_{ij})$, such that the rules $ {\bf W \to r, S \to q }$ 
make all the relations of the step 2, true.
Hereafter, we shall consider that all the polynomials, we work with in our proof, have been evaluated coherently to the above rules.
We have now, for the polynomials $A$ and $B$:
\[ A=\sum_{k=1}^{\omega} c_k(w_{ij}) \de_{{\bf i}_k} \ast L_k +R_{\de}=
\sum_{k=1}^{\omega} c_k(w_{ij}) \de_{{\bf i}_k} \ast \hat{L}_k \ast L +R_{\de}=
\hat{A} \ast L + R_{\de} \]

\[ B=\sum_{\rho =1}^{\omega} c_{\rho}(s_{ij}) \e_{{\bf j}_{\rho}} \ast M_{\rho} +R_{\e}=
\sum_{\rho =1}^{\omega} c_{\rho}(s_{ij}) \e_{{\bf j}_{\rho}} 
\ast \hat{M}_{\rho} 
\ast M +R_{\e}=\hat{B} \ast M +R_{\e}\]
Taking into account the fact that the
equalities of the step 2 are true, we can easily get that $\hat{A}=\hat{B}$,
$R_{\de}=R_{\e}=0$. 
The original system can now be re-written as 
$\hat{A} \ast Ly(t)=\hat{B} \ast M u(t)$. We define $y_1(t)=Ly(t)$,
$u_1(t)=Mu(t)$ and thus
$\hat{A}y_1(t)=\hat{B}u_1(t)$. Since $\hat{A}=\hat{B}$ and
$\hat{A},\hat{B}$ proper, we conclude that $y_1(t)=u_1(t)$, which means that $Ly(t)=Mu(t)$.
To finish the proof
we have to show that this linear system is causal too. In other words, we must prove that
$d(L)<d(M)$,
but from the causality
of the original system we get $d(A)<d(B)$ and $d(A)=d(\hat{A})+d(L)$, $d(B)=d(\hat{B})+d(M)$, $\hat{A}=
\hat{B}$ implies the desired result.
All the above statements hold on in the case of any other pair of rules 
$({\bf p,y})$, belonging to ${U}$ and thus, any set of polynomials
$L \in {\bf L}$, $M \in {\bf M}$ can form an exact linear-output linearization of the system
$Ay(t)=Bu(t)$.
The proof has been completed.
$\Box$
\begin{rem}\label{prin}
The main advantage of the above procedure is that it gives a set of solutions and it works 
when the polynomials $A,B$ have no linear parts.
Nevertheless, if the polynomials $A$ and $B$ do have linear parts, then two facts can be easily proved:
(a) the properness condition of step 2 is not needed
(b) the polynomials $L$ and $M$, are factors of the linear parts of $A$ and $B$.
\end{rem}
\begin{exa}
To clarify our ideas we present the following examples:
\par\noindent
(a) Let us consider the nonlinear discrete system:
\[ 4y(t)y(t-1)+2y(t)y(t-2)+22y^2(t-1)+21y(t-1)y(t-2)+5y^2(t-2)=\]
\[ = u(t-1)u(t-2)+2u^2(t-2)-3u(t-1)u(t-3)-21u(t-2)u(t-3)+45u^2(t-3) \]
or using $\de-operators$:
\[(4 \de_0\de_1+2\de_0\de_2+22\de_1^2+21\de_1\de_2+5\de_2^2)y(t)=\]
\[=(\e_1\e_2+2\e_2^2-3\e_1\e_3-21\e_2\e_3+45\e_3^2)u(t)\] shortly
$Ay(t)=Bu(t)$.
Applying the Linear-Equivalence Algorithm I, we take the following factorizations:
\[ A=5\de_0^2 \ast (w_{10}\de_0+w_{11}\de_1+\de_2)+(21-10w_{11})\de_0\de_1\ast(w_{20}\de_0+\de_1)+\]
\[+(22-5w_{11}^2-21w_{20}+10w_{11}w_{20})\de_0^2 \ast (w_{30}\de_0+\de_1) +R_{\de} \]
\[B=45\e_0^2 \ast (s_{10} \e_0 +s_{11} \e_1 +s_{12}\e_2+\e_3)+
(-21-90s_{12})\e_0\e_1 \ast (s_{20}\e_0+s_{21}\e_1+\e_2)+\]
\[+(-3-90s_{11}+21s_{21}+90s_{12}s_{21})\e_0\e_2 \ast (s_{30}\e_0+\e_1)+c_5'\e_0\e_1\ast(s_{50}\e_0+\e_1)+\]
\[+(2-45s_{12}^2+21s_{21}+90s_{12}s_{21})\e_0^2 \ast (s_{40}\e_0+s_{41}\e_1+\e_2)+c_6'\e_0^2
\ast (s_{60}\e_0+\e_1)+R_{\e}\]
(The terms $R_{\de},c_5',c_6',R_{\e}$ are not presented explicitly, due to their huge size.)
The set ${U}$, that is the set of the values of the parameters $w_{ij},s_{ij}$, which satisfy the
relations of Step 2, is the following:
${U}=\{(w_{10}=0,w_{11}=2,w_{20}=2,w_{30}=2,$
$s_{10}=0,s_{11}=0,s_{12}=-{1 \over 3}, s_{20}=0,s_{21}=-{1 \over 3},$ 
$s_{30}=0,s_{40}=0,s_{41}=-{1 \over 3},s_{50}=0,s_{60}=0),$
$(w_{10}=0,w_{11}=2,w_{20}=2,w_{30}=2,$
$s_{10}=0,s_{11}=0,s_{12}=-{1 \over 3}, s_{20}=0,s_{21}=-{1 \over 3},$ 
$s_{30}=-{1 \over 3},s_{40}=0,s_{41}=-{1 \over 3},s_{50}=-{1 \over 3},s_{60}=-{1 \over 3})\}$.
These values will give $\hat{A}=\hat{B}=5\de_1^2+\de_0\de_1$, which are proper polynomials, and
$L=2\de_0+\de_1$, $M=\e_1 -3 \e_2$.
Hence $U\ne \emptyset$ and the sets ${\bf L}$, ${\bf M}$ are 
${\bf L}=\{2\de_0+\de_1\}$, ${\bf M}=\{\e_1 -3 \e_2\}$. This means that the " linearization " of the above system is the linear system:
$(2\de_0+\de_1)y(t)=(\e_1 -3 \e_2)u(t)$ or
\[2y(t)+y(t-1)=u(t-1)-3u(t-2)\]
 To verify that the two systems, the nonlinear and the linear one, have the same dynamic 
 behavior we present, in table \ref{elin}, their responses under the same input sequence
 $u(t)=rnd(1)$, and the same initial conditions. 
  \par
  \begin{table}
\begin{center}
\begin{tabular}{||r|r@{.}l|r@{.}l||}
\hline
Repeats & \multicolumn{2}{|c|}{Non-linear Output} & \multicolumn{2}{|c|}{Linear Output} \\
\hline
5 & 0&158373 & 0&158373 \\
\hline
100 & -0&346706 & -0&346701  \\
\hline
300 & 0&024609 & 0&024609 \\
\hline
\end{tabular}
\end{center}
\caption{\label{elin} Simulation of the Simple Linear Simplification.}
\end{table}  
\noindent
(b) 
We have the causal nonlinear discrete system:
\[ 2y(t)+2y(t-1)+{1 \over 2} y(t-2)+y(t-1)y(t-2)-{1 \over 2}y^2(t-2)-{1 \over 2}y(t-1)y(t-3)-\]
\[ -{1 \over 4} y(t-2)y(t-3)+{1 \over 4}y^2(t-3)-{1 \over 2}y(t-1)y(t-4)+{1 \over 4}
y(t-2)y(t-4)+{1 \over 4}y(t-3)y(t-4)=\]
\[=2u(t-1)+u(t-2)+u(t-2)u(t-3)-u^2(t-3)-u(t-2)u(t-4)+u(t-3)u(t-4) \]
Using the same methodology, we can find that the linear system:
$(\de_0+{1\over 2} \de_1)y(t)=\de_1 u(t)$ or
\[y(t)+{1 \over 2}y(t-1)=u(t-1)\]
 consists the simplification upon request. We observe here that
the polynomials $A$ and $B$ have linear parts $A_l=2\de_0+2\de_1+\frac{1}{2}\de_2$, $B_l=2\e_1+\e_2$
and therefore the linear polynomials
$L=\de_0+{1\over 2} \de_1$, $M=\de_1$ are factors of $A_l$ and $B_l$ (Remark \ref{prin}).
\end{exa}
\subsection{The Cross-Products Case}
Let us suppose now that we have a causal non-linear 
discrete system with cross-products of the form
$Ay(t)=Bu(t)+{C}[y(t),u(t)]$. As we mentioned before, a linearization method for this system
consists in discovering linear systems that are equivalent to the original one.
The
algorithm for this category of systems is: 

\small {\sf \vskip 15 pt \noindent
\underline{\bf The Linear-Equivalence Algorithm II}
 \vskip 10 pt \noindent 
 {\bf Input:} The $\de$-polynomial $A$, the $\e$-polynomial $B$ and the $\de \e$-polynomial $C$. 
 \begin{quote}
{\bf STEP 0:} We define the $\de\e$-polynomial $\uA$ as
$\uA=A-B-C$.(We transfer all the terms at the left-hand side of (\ref{demain})).
\par
{\bf STEP 1:} By means of the $F\de \e L $-Algorithm, we find the Formal $\de \e
L$-Factorization of the polynomial $\uA$:
\[ \uA=\sum_{k=1}^{g} c_k(w_{ij},s_{ij})\de_{{\bf i}_k} \e_{{\bf j}_k} \ast [L_k,M_k] +R_{\uA}\]

\par
{\bf STEP 2:} We define the set of rules $U=\{(w_{ij},s_{ij})=$
$(a_{ij},b_{ij})$, $a_{ij},b_{ij} {\bf R}\}$ such
that the following are valid simultaneously
\begin{itemize}
\item
$\left. \begin{array} {c}
R_{\uA}\\
\end{array} \right|_{U}=0$, (we eliminate the remainder).

\item The sets $\bf{L}=\{{\it L:L}=\gcd\left( \left. \begin{array}
{c}
L_k\\
\end{array} \right|_{U}\right)\}$,
$\bf{M}=\{{\it{M}:M}=\gcd\left( \left. \begin{array} {c}
M_k\\
\end{array} \right|_{U}\right)\}$, $k=1,\ldots,g$,
are not void.
(We find the common factors of the linear polynomials).

\item The $\de \e$-polynomials, $\underline{\hat{A}}_\rho, \rho=1, \ldots, k$, where
$ \uA=(\sum_{\rho=1}^ka_\rho\underline{\hat{A}}_\rho) \ast [L,M]$, $L \in {\bf L}, M \in {\bf M}$, 
$a_\rho \in
{\bf R}$, are homogeneous, with respect to the  multiindices ${\bf i}_\rho, \rho=1, \ldots, k$ and
the standard constants $\lambda, \mu$. (We factorize $\uA$ with respect
to the common factors).

\item The polynomial 
$\hat{A}=\sum_{\rho=1}^k a_\rho \de_{{\bf i}_\rho}$,
is proper.
(The multiindices ${\bf i}_{\rho}$ are the multiindices of the homogeneous polynomial).

\end{itemize}
\par
{\bf STEP 3}: {\bf IF} $U \ne \emptyset$ {\bf THEN} goto the output {\bf ELSE} the
method fails.

\end{quote}
} 
\vskip 15 pt
{\bf Output:} The quantities
$\bf{L}, \bf{M}$, $\lambda, \mu$

\normalsize
\par
\noindent
The condition $U \ne \emptyset$ is equivalent to the solvability of an
algebraic system
of polynomial equations. Indeed, since the remainder
$R_{\uA}$ is a $\de\e$-polynomials with coefficients which involve 
the parameters $(w_{ij},s_{ij})$,
the condition $\left. \begin{array} {c}
R_{\uA}\\
\end{array} \right|_{U}=0$, is equivalent to a system of polynomial equations.
This system can be solved by successive substitutions. The first equation contains only one unknown parameter and thus it is solvable. Substituting the value will obtain to the second equation we take a solution for the second parameter and so on. Among the solutions of this algebraic system, if any, we can choose those 
particular one, which can satisfy and other conditions, like minimality
in the number of evaluated parameters and so on.

\par
\begin{theo}
Let $Ay(t)=Bu(t)+{C}[y(t),u(t)]$, $A$ a $\de$-polynomial, $B$ an $\e$-polynomial and
$C$ a $\de \e$ - polynomial, be a causal non-linear
discrete system that contains products among input and outputs. If ${\bf
L,M}$,$\lambda$, $\mu$, are the outputs of the previous
algorithm, then, any linear system of the form: $\lambda Ly(t)=-
\mu Mu(t)$, $L \in {\bf L}, M \in {\bf M}$ is
equivalent to the original non-linear system.
\end{theo}
{\bf Proof:} Let ${\bf W}=\{w_{ij}\}$ and ${\bf S}=\{s_{ij}\}$, be the sets of the parameters 
 that appear in the
Formal $\de \e L$-Factorization of $\uA$, during the run of 
the algorithm. The existence
of the solution entails that ${U} \ne \emptyset$ and, therefore, 
that there is at
least one pair of vectors, ${\bf r}=( a_{ij})$, ${\bf q}=( b_{ij})$, such that the
rules $ {\bf W \to r, S \to q }$ make all the relations of the step 2, true.
Hereafter, we will assume that all the polynomials in the proof, have been evaluated
and cohere with the above rules. For the polynomial $\uA$ we now have:
\[ \uA=\sum_{k=1}^{g} c_k(w_{ij},s_{ij})\de_{{\bf i}_k} \e_{{\bf j}_k} \ast [L_k,M_k] +R_{\uA}=
\sum_{k=1}^{g} c_k(w_{ij},s_{ij})\de_{{\bf i}_k} \e_{{\bf j}_k} \ast [\hat{L}_k\ast L,\hat{M}_k \ast M] +R_{\uA}=\]
\[=(\sum_{\rho=1}^ka_\rho\underline{\hat{A}}_\rho) \ast [L,M]+R_{\uA}\]
By taking into account the fact that the equalities
in step 2 are true, 
we can easily see that $R_{\uA}=0$. 
Moreover the polynomials $\underline{\hat{A}}_\rho$ are homogeneous.
This means that taking into account proposition \ref{roundpro} we get
$\uA=(\sum_{\rho=1}^k a_\rho \de_{{\bf i}_\rho}) \ast (\lambda L+\mu M)$, 
or $\uA=\hat{A} \ast (\lambda L+\mu M)$.
The original system can be now
re-written as
$\hat{A} \ast (\lambda L+\mu M)[y(t),u(t)]=0$
We define
$ \omega(t)=\lambda Ly(t)+\mu Mu(t)$,
and the system becomes
$\hat{A} \omega(t)=0$. The last condition
of the step 2, guarantee that $\hat{A}$ is proper.
This means that $\omega(t)=0$ (corollary \ref{porisma} ) or
$\lambda Ly(t)+\mu Mu(t)=0$ and thus
$\lambda L y(t)=-\mu Mu(t)$. To conclude the proof we need to show that this linear
system is also causal. In other words, we must prove that 
$d(\lambda L)<d(-\mu M)$.
From the causality of the original system we get:
$d(\uA_y)<d(\uA_u)$, where
$d(\uA_y)$ is the minimum delay of the output signal at the polynomial $\uA$,
$d(\uA_u)$ is the minimum delay of the input. 
Additionally, using the equalities, presented above, we get:
$d(\uA_y)=d(\hat{A})+d(L)$, 
$d(\uA_u)=d(\hat{A})+d(M)$.
The elimination of $d(\hat{A})$
entails
the desired result. All the above
statements hold for any other pair of rules $({\bf p,y})$ belonging to ${U}$ and thus,
any set of polynomials $L \in {\bf L}$, $M \in {\bf M}$ can form an exact
linear-output linearization of the nonlinear system. This completes the proof.
$\Box$
\begin{exa}
We have the system: 
\[ 6y(t)-11y(t-1)+6y(t-2)-y(t-3)+36y(t-1)y(t-2)-\]
\[-30y(t-1)y(t-3)-30y^2(t-2)+25y(t-2)y(t-3)-y(t-3)y(t-4)=\]
\[=u(t-1)+u(t-2)-2u(t-3)-u(t-2)u(t-3)+4u(t-3)u(t-4)+\]
\[+6y(t-1)u(t-3)-5y(t-2)u(t-3)+6u(t-2)y(t-2)-\]
\[-5y(t-3)u(t-2)-2u(t-3)y(t-4)-2y(t-3)u(t-4)\]
 Putting all the terms at the left-hand side of the above equation
 and  using $\de$,$\e$ and $\de\e$-operators we find $\uA$. That is: 
 \[ \uA=6\de_0-11\de_1+6\de_2-\de_3-\e_1-\e_2+2\e_3+36\de_1\de_2-30\de_1\de_3-\]
\[-6\de_1\e_3-30\de_2^2+25\de_2\de_3+5\de_2\e_3-6\de_2\e_2+5\de_3\e_2+\]
\[+\e_2\e_3-4\e_3\e_4+2\e_3\de_4+2\de_3\e_4-\de_3\de_4\]
Following the Linear-Equivalence II Algorithm
step by step we get $U \ne \emptyset$ for a specific set of values
for the parameters, denoted as $({\bf r,q})$. (For reasons of brevity, we have not shown all
calculations in detail). Using the substitution rule $({\bf W,S}) \to ({\bf r,q}) $,
we rewrite the polynomial $\uA$ as follows:
\[ \uA=\de_0\e_e \ast [6\de_0-5\de_1+\de_2,-\e_1-2\e_2]+\de\e_0 \ast[6\de_0-5\de_1+\de_2,-\e_1-2\e_2]-\]
\[-\de_1\e_e \ast[6\de_0-5\de_1+\de_2,-\e_1-2\e_2]
 -\de_e \e_1 \ast[6\de_0-5\de_1+\de_2,-\e_1-2\e_2]+\]
\[+\de_1\de_2 \e_e \ast[6\de_0-5\de_1+\de_2,-\e_1-2\e_2]
+\de_1\e_2 \ast[6\de_0-5\de_1+\de_2,-\e_1-2\e_2]+\]
\[+\de_2\e_1 \ast[6\de_0-5\de_1+\de_2,-\e_1-2\e_2]+\de_e\e_1\e_2 \ast[6\de_0-5\de_1+\de_2,-\e_1-2\e_2]=\]
\[=[(\de_0\e_e+\de_e\e_0)-(\de_1\e_e+\de_e\e_1)+(\de_1\de_2\e_e+\de_1\e_2+\de_2\e_1+\]
\[+\de_e\e_1\e_2)] \ast[6\de_0-5\de_1+\de_2,-\e_1-2\e_2]\] The polynomials in the parenthesis on the right of the 
 above expression, are homogeneous  
polynomials in relation to the multiindices, 
$0,1,(1,2)$ and constants $\lambda =\mu =1$.
We see that $a_0=b_0=1$, $a_1=b_1=-1$, $a_{(1,2)}=b_{(1,2)}=1$
and therefore we construct the polynomial
$\hat{A}=\de_0-\de_1+\de_1\de_2$, which is proper.
Thus, finally we get ${\bf L}=\{ 6
\de_0-5\de_1 \}$, ${\bf M}=\{ -\e_1 \}$, $\lambda =\mu =1$. Therefore the desired linear system is
$Ly(t)=-Mu(t)$ or $(6\de_0-5\de_1+\de_2)y(t)=(2\e_2+\e_1)u(t)$ or 
\[6y(t)-5y(t-1)+y(t-2)=u(t-1)+2u(t-2)\]
 Simulations are presented in Table \ref{mixlin}.
We took $u(t)=rnd(1)$ as input.
\begin{table}
\begin{center}
\begin{tabular}{||r|r@{.}l|r@{.}l||}
\hline
Repeats & \multicolumn{2}{|c|}{Non-linear Output} & \multicolumn{2}{|c|}{Linear Output} \\
\hline
50 & 0&922726 & 0&922726 \\
\hline
500 & 0&656207& 0&656207 \\
\hline
1000 & 1&15928 & 1&15928 \\
\hline
\end{tabular}
\end{center}
\caption{\label{mixlin} Simulation of the Linear Equivalence}
\end{table}
\noindent

\end{exa}

\section{Concluding Remark}

The aim of this paper was to describe algebraic computational methods for the simplification of
a general class of non-linear discrete input - output systems that contain products
between input and output signals. Actually, we developed an approach to the
Linear - Equivalence problem. The entire approach is based on a
proper framework that involves the so-called $\de\e$-operators and the star-product
operation. We hope to be able to present current work on further applications of this
method to more concrete questions in a future paper.

\end{document}